\documentclass
[manmat
]
{svjour}

\usepackage[dvips]{graphicx}

\usepackage{amsmath}
\usepackage{amssymb}

 \spnewtheorem*{HT}{Hamburger Theorem}{\bf}{\it}

 \spnewtheorem*{Conj}{Erd\"os
Conjecture}{\bf}{\it}

\def\Zer#1{\mathcal{Z}(#1)}
\def\Com#1{\mathbb{C}^{#1}}

\def\R#1{\mathbb{R}^{#1}}
\def\Z#1{\mathbb{Z}^{#1}}
\def\Tcirc{\mathbb{T}}
\def\Uu{\mathbb{D}}
\def\Uv{\mathbb{U}}
\def\Haus#1{\mathcal{H}^{#1}}

\def\Del{\frak{I}}

\def\Dd{\mathcal{U}}

\def\z{\zeta}

\def\scal#1#2{\langle #1; #2 \rangle}

\def\vv{\vec{v}}
\def\conj#1{\overline{#1}}

\def\scalm#1#2{\scal{#1}{#2}_{t}}

\def\Dp{\mathcal{D}}

\DeclareMathOperator{\Div}{div} 
 \DeclareMathOperator{\re}{Re}

\DeclareMathOperator{\length}{length}

\begin{document}


\title{Length functions of lemniscates}

\author{
\inst{1}Kuznetsova O.S., \thanks{\inst{1}The author was supported
the G\"oran Gustafsson foundation and grant RFBR no.~03-01-00304.}
\and \inst{2}Tkachev V.G.
\thanks{\inst{2}The author was supported by Russian President grant
for young doctorates no.~ 00-15-99274 and grant RFBR
no.~03-01-00304.} }

\authorrunning{Kuznetsova O.S. and Tkachev V.G.}

\institute{Kuznetsova O.S., Tkachev V.G. \at Royal Institute of
Technology, Lindstedtsv\"agen 25, 10044 Stockholm, Sweden, \and
Volgograd State University, 2-ya Prodolnaya 30, Volgograd, 400062
Russia \\(\email{astra1987@mail.ru},
\email{tkatchev@math.kth.se})}

%
%

%
\subclass{30E05, 42A82, 44A10}

\date{Received: date / Revised version: date }

\addtolength{\baselineskip}{0.4mm}



\maketitle

\begin{abstract}
We study metric and analytic properties of generalized lemniscates
$E_t(f)=\{z\in\Com{}:\ln |f(z)|=t\}$, where $f$ is an analytic
function. Our main result states that the length function
$|E_t(f)|$ is a bilateral Laplace transform of a certain positive
measure.  In particular, the function $\ln |E_t(f)|$ is convex on
any interval free of critical points of $\ln|f|$. As another
application we deduce explicit formulae of the length function in
some special cases.
\end{abstract}

\keywords{Lemniscate,  exponentially convex function, completely
monotonic function, Hamburger moment problem, harmonic level sets}



\section{Introduction}
Throughout this paper $E_f(t)$ denotes the $t$-level set
\begin{equation}\label{lemnisc}
\ln |f(z)|=t
\end{equation}
of an analytic function $f(z)$.

Our starting point is the polynomial lemniscates. Let $f(z)$ be a
monic polynomial $P(z)=z^n+a_1z^{n-1}+\ldots+a_n$, $n\geq 2$. In
1958 Erd\"os, Herzog and Piranian \cite{EHP} posed a number of
problems concentrated around the metric properties of lemniscates
(see also the later paper \cite{Erd}). Among them is the following

\begin{Conj}[Problem~ 12, \cite{EHP}; Problem VI, \cite{Erd}]
For fixed degree $n$ of $P$, is the length of the lemniscate
$|P(z)|=1$ greatest in the case where $P(z)=Q_n(z):=z^n-1$? Is the
length at least $2\pi$, if $E_P(0)$ is connected?
\end{Conj}

The actual breakthrough in the Erd\"os conjecture was made
recently by A.~ Eremenko and W.~ Hayman in \cite{ErHay}. They
proved that for any degree $n$ there exist a polynomial $P^*(z)$
which maximizes the length of $E_P(0)$ (following \cite{ErHay}, we
call it an \textit{extremal} polynomial); moreover, the estimate
\begin{equation*}\label{er1}
c_n\equiv \max_{\deg P=n}|E_P(0)|\leq A n
\end{equation*}
holds, where $|E|$ denotes the length of $E$, and $A\approx
9.173$. One can readily check that the conjectural value is $
|E_{Q_n}(0)|=2n+O(1)$. The previous upper estimates were due to
Ch.~Pommerenke \cite{Pom61}: $c_n\leq 74 n^2$ and P.~Borwein
\cite{Bor}: $c_n\leq 8\pi e n$.

Another important results of \cite{ErHay} states that \textit{the
lemniscate $E_{P^*}(0)$ is always connected and for any degree $n$
there exists an extremal polynomial $P^*$  such that all its
critical points belong to $E_{P^*}(0)$.}

Concerning the first part of Erd\"os conjecture, which is still
unsolved, T.~ Erdelyi writes that ``this problem seems almost
impossible to settle'' \cite[p.~8]{Erdelyi}. Another difficulty in
the study of the problem is the absence of any explicit formulae
for the \textit{length function} $|E_P(t)|$ (except for the
trivial polynomials $P(z)=(z-a)^n$). This question was initially
posed by Piranian in \cite{Piran} for the rose-type polynomials
$Q_n$; explicit formulae of $|E_{Q_n}(t)|$ were obtained by Butler
\cite{Butler} and Elia \cite{Elia}.

The second part of Erd\"os conjecture is related to the lower
estimate of $|E_P(0)|$ for so-called $K$-polynomials, i.e. the
polynomials with connected lemniscate $E_P(0)$. This problem was
solved in affirmative by Pommerenke in \cite{Pom59-2}, who
established that
\begin{equation}
\label{er0} \min_{ P\in K,\; \deg
P=n}|E_P(0)|=|E_{(z-a)^n}(0)|=2\pi.
\end{equation}

\subsection{Main results}\label{subsec:intro}

\begin{figure}[ht]
\includegraphics[height=0.25\textheight,keepaspectratio=true]{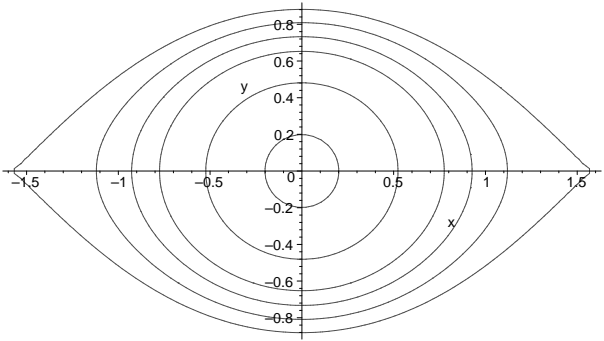}
\caption{The regular lemniscate region for $f(z)=\sin z$}
\label{fig:1}\end{figure}

\begin{definition}
By a \textit{lemniscate region} of $f$ we mean a triple
$(\Dd{},f,\Del)$ where $f(z)$ is an analytic function, $\Dd{}$ is
a component of the set
$$
\Dd{}_f(\Del):=\{z\in\Com{}:a<\ln|f(z)|<b\},
$$
$\Del=(a,b)$, such that for every $t\in \Del$ the set
$$
E_{\Dd{},f}(t):=E_f(t)\cap \Dd{}
$$
is compact in $E_f(t)$. A lemniscate region will called
\textit{regular} if $\Dd{}$ contains no zeroes and no critical
points, i.e.$f(z)f'(z)\ne 0$ in $\Dd{}$ (cf.
\cite[p.~264]{Hille}).
\end{definition}
It is easy to see that any (analytic) function $f$ admits regular
lemniscate regions provided that  the set of zeroes
$$
\Zer{f}:=\{z:\;f(z)=0\}
$$
is nonempty. Disregarding the polynomial case, where all
lemniscates are compact curves, we mention, for instance, an
example of a regular lemniscate region (the corresponding
lemniscate family is drawn on Figure~\ref{fig:1})
$$
f(z)=\sin z, \qquad \Del=(-\infty,0), \qquad \Dd{}=\{z: 0<|\sin
z|<1, |\re z|<\frac{\pi}{2}\}.
$$

By its definition, the projection $\ln|f(z)|:\Dd{}\to \Del$ is a
proper map. Hence, in case of a regular lemniscate region, the
level set $E_{\Dd{},f}(t)$  splits into a finite collection of
simple closed curves; any finite collection of closed components
of $E_{\Dd{},f}(t)$ will be called a $t$-\textit{lemniscate}, or
just a lemniscate of $f$.

Given a lemniscate domain $(\Dd{},f,\Del)$ we  define the
corresponding \textit{length function}
\begin{equation}\label{length} |E_{\Dd{},f}(t)|:=\length
(E_{\Dd{},f}(t)), \quad t\in\Del{},
\end{equation}
which will be in focus of the present paper. On the other hand,
the expressions for the higher derivatives of $|E_{\Dd{},f}(t)|$
(see (\ref{equ:H-prime}) below) involve integrals of the kind
\begin{equation}\label{def-s}
L_w(t)=\|w\|^2_t:=\int\limits_{E_{\Dd{},f}(t)}|w(z)|^2\; |dz|,
\end{equation}
where $w(z)$ is an analytic in $\Dd{}$ function. This makes it
natural to consider the averages (\ref{def-s}) as a suitable
generalization of the length function. The case $w\equiv 1$
obviously reduces to the length function.

Further, we consider the following first-order differential
operator
\begin{equation}\label{equ:w-def}
G_f(w)=w_{[1]}:=2gw'+g'w, \qquad w_{[k]}=G_f^{k}(w),
\end{equation}
where $g=f/f'$, and $w_{[k]}$ are the $G_f$-iterations of $w\equiv
w_{[0]}$.

The following theorem summarizes the results given in
Section~\ref{sec:proof} below.

\begin{theorem}\label{theo:H-main}
Let $(\Dd{},f,\Del)$ be a regular lemniscate region and $w(z)$ be
an analytic in $\Dd{}$ function. Then $L_w(t)\in
C^{\infty}(\Del)$. Moreover, for any $t\in\Del$ the sequence of
 derivatives $L_w^{(k)}(t)$, $k\geq 0$, $t\in\Del$, forms
a Hamburger moment sequence; that is, for any $t\in\Del$ and for
all $p=0,1,\ldots$, the Hankel matrices
\begin{equation}\label{HAN}
\left(%
\begin{array}{cccccc}
L_w(t) & L_w'(t) & \ldots & L_w^{(p)}(t)\\
L_w'(t) & L_w''(t)  &  \ldots &L_w^{(p+1)}(t) \\
\vdots & \vdots & \ddots & \vdots \\
L_w^{(p)}(t) & L_w^{(p+1)}(t)&\ldots  &L_w^{(2p)}(t)
\end{array}%
\right)
\end{equation}
are non-negative definite. Moreover, the entries of the latter
matrix are the following scalar products
\begin{equation}\label{equ:s-k-derivative}
L_w^{(k)}(t)=\scal{w_{[j]}}{{w}_{[k-j]}}_t:=\int\limits_{E_{\Dd{},f}(t)}\re
\overline{w}_{[k-j]}w_{[j]}\;|dz|, \qquad 0\leq j \leq k
\end{equation}
where the  integrals are independent of the choice of $j$.
\end{theorem}

Now, applying the well-known Bernstein's theorem \cite{Bern-AMF}
(see, also \cite[Ch.~ VI]{Widder}), we obtain a bilateral Laplace
representation of $L_w$

\begin{corollary}\label{theo:B-W-rep}
Let $(\Dd{},f,\Del)$ be a regular lemniscate region. Given an
analytic in $\Dd{}$ function $w(z)$ there exists a non-decreasing
function $\sigma(x)$  such that
\begin{equation}\label{equ:Bern-Wid}
L_w(t)=\int\limits_{E_{\Dd{},f}(t)}|w(z)|^2\;
|dz|=\int\limits_{-\infty}^{+\infty}e^{xt}\; d\sigma(x), \qquad
t\in (\alpha;\beta),
\end{equation}
and the latter integral converges for $t\in\Del$.
\end{corollary}

Functions $L(t)$ which satisfy a bilateral Laplace representation
(\ref{equ:Bern-Wid}) are known as \textit{exponentially convex}
functions (\cite{Bern-AMF}, \cite[\S~V.5.4]{Akhiezer}). This means
that the associated with $L(t)$ stationary kernel
$\mathbf{L}(x,y)=L(\frac{x+y}{2})$ is of positive type, i.e. for
every finite sequence $\{t_j\}_{1}^{m}$ from $\Del$ the quadratic
form
$$
\sum_{i,j=1}^{m} L\left(\frac{t_i+t_j}{2}\right)\xi_i\xi_j
$$
is positive (definite or semidefinite). In particular, given an
exponentially convex function $L(t)$ the function $\ln L(t)$ is
convex.

This class was introduced and extensively studied  by S.~Bernstein
\cite{Bern-AMF} and D.~Widder \cite{Widder} in connection with the
so-called completely (or absolutely) monotonic analytic functions
(see the definition in Section~\ref{sec:measure}). We only mention
a deep penetration of the both classes into complex analysis,
inequalities analysis \cite{Alzer}, special functions
\cite{Samko}, probability theory \cite{Kimb}, radial-function
interpolation \cite{Umem}, harmonic analysis on semigroups
\cite{semi} (for further discussion and references, see recent
survey \cite{BergDuran}).

We also remark that exponential convexity leads to further
 inequalities on $L_w$  and its derivatives like
those considered in \cite{GPP}.

As another consequence of exponential convexity (see
\cite[\S~15]{Bern-AMF}) we point out the following continuation
property, actually, a \textit{complexification} of the length
function (see also (\ref{complexific}) below)

\begin{corollary}\label{corol:cont}
Under the hypotheses of Theorem~\ref{theo:H-main}, the function
$L_w(t)$ admits an analytic continuation $L_w(z)$ into the strip
$a<\re z<b$, where $\Del=(a,b)$.
\end{corollary}

In the remaining part of this section we consider the case when
the function $f(z)$ is  a monic polynomial $P(z)$ and $w(z)\equiv
1$. By $T_1<\ldots <T_{\nu-1}$ we denote the set of all (finite)
critical values of $\ln |P|$. The intervals
$\Del_j=(T_{j-1},T_{j})$ will be called \textit{regular} (with
respect to $P$), where $T_0=-\infty$, $T_\nu=+\infty$. Then
$(\Dd_P(\Del_j),P,\Del_j)$, $1\leq j\leq \nu-1$, constitutes a
special class of the \textit{principal} regular lemniscate
regions.

We have in the preceding notations

\begin{corollary}\label{corol:exp-conv}
Given a regular interval $\Del_j=(T_{j-1},T_{j})$, $1\leq j\leq
\nu$, the following representation holds
\begin{equation}\label{equ:rr}
|E_P(t)|=\int\limits_{-\infty}^{+\infty}e^{xt}\;
d\sigma^{P,\Del_j}(x), \qquad t\in\Del_j,
\end{equation}
where $\sigma^{P,\Del_j}(x)$ is a non-decreasing function. In
particular, $\ln|E_P(t)|$ is a convex function on $\Del_j$.
\end{corollary}

It turns out (see Section~\ref{sec:measure} below)  that the
function $\sigma^{P,\Del_\nu}(x)$ (i.e., $j=\nu$) is a piece-wise
constant function and the integral in (\ref{equ:rr}) can be
written as sum of a certain exponential series.

Given a monic polynomial $P$, $\deg P=n$, we define an auxiliary
function
\begin{equation*}\label{Phi} \Phi_P(t):=\ln
|E_P(t)|-\frac{t}{n}.
\end{equation*}
We call $\Phi_P(t)$ the \textit{indicator} of $P$. Then
$\Phi_P(t)$ has a simple invariance  property with  respect to
dilatations of $P$
\begin{equation}\label{equ:transf}
 P_\alpha(z):=
e^{-\alpha}P(ze^{\frac{\alpha}{n}})=z^n+a_1e^{-\frac{\alpha}{n}}z^{n-1}+\ldots
+a_{n-1}e^{-\frac{\alpha(n-1)}{n}}z+a_{n}e^{-\alpha}
\end{equation}
where $\alpha\in (-\infty;+\infty]$.  Indeed, we note that
$P_\alpha(z)$ is also a monic polynomial, $P_0(z)\equiv P(z)$, and
\begin{equation}\label{equ:homot}
E_{P_\alpha}(\beta)=e^{-\frac{\alpha}{n}}E_P({\alpha+\beta}).
\end{equation}
Hence, we have
\begin{equation}\label{equ:Phi-homot}
\Phi_{P_\alpha}(\beta)=\Phi_{P}(\alpha+\beta).
\end{equation}

\begin{proposition}\label{rem:1} Let $P^*(z)$ be an extremal
(with respect to Erd\"os conjecture) polynomial of degree $n$;
then $t=0$ is an absolute \textit{maximum} point of
$\Phi_{P^*}(t)$.
\end{proposition}

\begin{proof}
It follows from (\ref{equ:Phi-homot}) that
$$
\Phi_{P^*}(t)=\Phi_{P^*_t}(0)=|E_{P^*_t}(0)|\leq
|E_{P^*}(0)|=\Phi_{P^*}(0)
$$
which proves the required property.
\qed\end{proof}

\begin{theorem}\label{th-er1}
$\Phi_P(t)$ and $|E_P(t)|$ are continuous functions in
$(-\infty;+\infty)$. Moreover, if the polynomial $P(z)$ is
non-trivial \textrm{(}i.e. is different from $(z-a)^n$\textrm{)}
then
\begin{enumerate}
    \item[A)] $\Phi_P(t)$ and $|E_P(t)|$ are strictly convex in
    each regular interval $(T_{j-1},T_j)$, $1\leq j\leq \nu$;
    \item[B)] the following asymptotics behavior holds
\begin{equation}\label{limits}
\lim_{t\to+\infty}\Phi_P(t)=\ln 2\pi.
\end{equation}
\end{enumerate}
\end{theorem}

\begin{remark}
A direct analysis near critical points of $\ln|P|$ implies  that
$\Phi_P(t)$ is only of H\"older class there which makes useless
the standard variational methods near the corresponding extremum.
\end{remark}

As another  consequences of Theorem~\ref{th-er1} we mention the
 estimate (\ref{er0}) due to Pommerenke.

\begin{corollary}\label{corol:ex}
Let $P$ be a monic $K$-polynomial, i.e. the lemniscate $E_P(0)$ is
connected. Then
$$
|E_P(0)|\geq 2\pi,
$$
with equality only in the case $P(z)=(z-a)^n$.
\end{corollary}

\begin{proof}
We can assume that $P(z)\ne(z-a)^n$, otherwise $|E_P(0)|=2\pi$.
Then it easily follows from the definition of $K$-polynomial that
$T_k\leq 0$ for all $k\leq \nu$. Hence, by virtue of
Theorem~\ref{th-er1} we conclude that $\Phi_P(t)$ is strictly
convex in $[0;+\infty)$. Moreover, by (\ref{limits}) the function
$\Phi_P(t)$ is bounded on $[0;+\infty)$ and it follows from the
strict convexity of $\Phi_P(t)$ that it is actually strictly
decreasing. Because of $\Phi_P(0)=\ln |E_P(0)|$, we have
$\Phi_P(0)>\ln 2\pi$, or $|E_P(0)|>2\pi$, which completes the
proof. \qed\end{proof}

By the  Eremenko-Hayman theorem, we know that for all integers
$n\geq 2$ an extremal polynomial $P^*$, $\deg P=n$, does exist
such that $\Zer{ {P^*}'}\subset E_P(0)$. The following assertion
gives a complement to the latter property.

\begin{corollary}\label{corol:level}
If $P^*(z)$ is an extremal polynomial of degree $n$ then the
lemniscate $E_{P^*}(0)$ is singular, i.e. it contains at least one
critical point:
$$
\Zer{ {P^*}'}\cap E_{P^*}(0)\ne \varnothing.
$$
\end{corollary}

\begin{proof}
By Proposition~\ref{rem:1}, $t=0$ is an absolute maximum of the
indicator function $\Phi_{P^*}(t)$ and it follows that $t=0$ can
not be a regular value of $\ln |P(z)|$ for $\Phi_{P^*}(t)$ is
strictly convex in a neighborhood of regular values.
\qed\end{proof}

As another application of (\ref{equ:s-k-derivative}), in
Section~\ref{subsec:examples} we obtain explicit formulae for the
length functions in the case when $f$ is a solution of the
following equation
\begin{equation*}\label{equ:generalized}
\varphi '=C(1-\varphi ^\nu)^{\frac{k+1}{\nu}}.
\end{equation*}

\section{Proof  of the main results}\label{sec}

\subsection{Preliminaries} \label{sec:aver}
Here we prove the main technical result which we formulate in a
form suitable for further applications.

Let $M$ be a $p$-dimensional Riemannian manifold and by
$\scal{X}{Y}$ and $\nabla$ the intrinsic scalar product and
covariant derivative are denoted. By $\Div X$ we denote the
divergence of a vector field $X$ generated by $\nabla$. We recall,
that a function $u(x):M\to \R{}$ is called \textit{harmonic} if
$\Delta u\equiv \Div \nabla u(x)=0$; by $\Sigma_u(t)$ we denote
the level set $\{x\in M: u(x)=t\}$.

\begin{definition}
A triple $(\Dd{},u,\Del  )$, where $\Dd{}$ is an open subset,
$\overline{\Dd{}}\subset M$, $u(x)$ is a harmonic function in
$\Dd{}$, and $\Del =(\alpha;\beta)$, is said to be a
\textrm{{(regular) lemniscate region}} if for all $t\in\Del{}$ the
set $\Sigma_{\Dd{},u}(t):=\Sigma_u(t)\cap \Dd{}$ is compactly
contained in $\Dd{}$ and $\Dd{}$ is free of critical points of
$u$.
\end{definition}

\begin{remark}\label{rem:nota}
Clearly, when $M=\Com{}$ the preceding definition is reduces to
the that one given in Section~\ref{subsec:intro}. We will not
distinguish the corresponding notations $(\Dd{},f,\Del )$ and
$(\Dd{},\ln|f|,\Del )$, when $f$ is an analytic function;
moreover, in this case
$$
E_{f}(t)=\Sigma_{\ln|f|}(t), \qquad
E_{\Dd{},f}(t)=\Sigma_{\Dd{},\ln|f|}(t).
$$
\end{remark}

\begin{lemma}\label{theo:exhaust}
Let $(\Dd{},u,\Del)$ be a lemniscate region and  $h(x)$ be a
$C^2$-smooth function on $\Dd{}$ and
\begin{equation}\label{equ:ex1}
H(t):=\int\limits_{\Sigma_{\Dd{},u}(t)}h(x)|\nabla u(x)|\,
d\Haus{p-1}(x),
\end{equation}
where $d\Haus{p-1}$ denotes $(p-1)$-dimensional Hausdorff measure
on ${\Sigma_{\Dd{},u}(t)}$.

Then $H(t)\in C^2(\Del) $ and
\begin{equation}\label{equ:ex2}
H'(\tau)=\int\limits_{\Sigma_{\Dd{},u}(t)}\frac{\scal{\nabla
h(x)}{\nabla u(x)}}{|\nabla u(x)|}\, d\Haus{p-1}(x),
\end{equation}
\begin{equation}\label{equ:ex3}
H''(\tau)=\int\limits_{\Sigma_{\Dd{},u}(t)}\frac{\Delta
h(x)}{|\nabla u(x)|}\, d\Haus{p-1}(x).
\end{equation}
\end{lemma}

\begin{proof}
Because of regularity condition  all the level sets
$\Sigma_{\Dd{},u}(t)$ are embedded submanifolds in $M$ and the
vector field
\begin{equation}\label{equ:nu}
\vec{\nu}(x)\equiv \frac{\nabla u(x)}{|\nabla u(x)|}
\end{equation}
represents the field of unit normals to  $\Sigma_{\Dd{},u}(t)$
pointing in the growth direction of $u$,
\begin{equation}\label{equ:normal}
\scal{\vec{\nu}(x)}{\nabla u(x)}=|\nabla u(x)|.
\end{equation}

We claim that for any $C^1$-vector field $\vv$ on $\Dd{}$
\begin{equation}\label{equ:Ftau}
\frac{d}{d\tau}\int\limits_{\Sigma_{\Dd{},u}(\tau)}\scal{\vv}{\vec{\nu}}\,d\Haus{p-1}=
\int\limits_{\Sigma_{\Dd{},u}(\tau)}\frac{\Div \vv}{|\nabla
u|}\,d\Haus{p-1}.
\end{equation}

Indeed, let $t\in \Del$, $t\ne \tau$, be chosen arbitrary. Then by
virtue of (\ref{equ:normal}) and harmonicity of $u(x)$ we have by
Stokes' formula
\begin{equation}\label{equ:Stokes}
\begin{split}
F(t)-F(\tau)&=\int\limits_{\Dd{}(t)-\Dd{}(\tau)}\scal{\vv}{\vec{\nu}}\,
d\Haus{p-1}=\\
&=\int\limits_{\partial \Dd{}(\tau,t)}\scal{\vv}{\vec{\nu}}\,
d\Haus{p-1}=\int\limits_{\Dd{}(t)-\Dd{}(\tau)}\Div \vv \;d x,
\end{split}
\end{equation}
where $\Dd{}(t)=\{x\in\Dd{}:u(x)<t\}$ and $F(\tau)$ denotes the
left-hand side integral in (\ref{equ:Ftau}).

Then applying co-area formula
to (\ref{equ:Stokes}) we obtain
\begin{equation}\label{equ:Leibniz}
\frac{F(t)-F(\tau)}{t-\tau}=\frac{1}{t-\tau}
\int\limits_{\tau}^{t}d\xi\int\limits_{\Sigma_{\Dd{},u}(\xi)}\frac{\Div
\vv}{|\nabla u|} \,d\Haus{p-1}.
\end{equation}

The latter limit does exist for every regular value $\tau$ of
$u(x)$ (even if $u$ is only locally Lipschitz in $D$
\cite[\S~3.2]{Federer}) and  (\ref{equ:Ftau}) follows.

Thus, applying (\ref{equ:Ftau}) to $\vv=h\nabla u$ we obtain
(\ref{equ:ex2}). Moreover, it follows from (\ref{equ:nu}) that
(\ref{equ:ex2}) can be written in the form
$$
H'(\tau)=\int\limits_{\Sigma_{\Dd{},u}(\tau)}\scal{\nabla
h}{\vec{\nu}}\, d\Haus{p-1}
$$
and, again applying (\ref{equ:Ftau}) to the last relation, now
with $\vv=\nabla h(x)$ we arrive at
\begin{equation*}
H''(\tau)=\int\limits_{\Sigma_{\Dd{},u}(\tau)}\frac{\Div \nabla
h}{|\nabla u|}\, d\Haus{p-1}
\end{equation*}
 and the lemma is proved.
\qed\end{proof}

\begin{remark}\label{remmm}
We notice that in the case when $M$ is a minimal submanifold of
$\R{N}$ and $u(x)$ is a coordinate function on $M$, the definition
of a regular lemniscate region corresponds to a special class of
minimal surfaces, so-called minimal tubes, in Euclidean space
\cite{Mik79}, \cite{MT89}. In that case a result similar to
Lemma~\ref{theo:exhaust} (actually, for the radial symmetric
functions $h$) was obtained by V.~Klyachin in \cite{Klya}.
\end{remark}

The following assertion is an easy consequence of Cauchy's
inequality and the last theorem, and it can be regarded as a
particular case of Theorem~\ref{theo:H-main} for $p=1$.

\begin{corollary} Let $\Delta \ln h\geq 0$; then, under the hypotheses of
Lemma~\ref{theo:exhaust},
 the function $\ln H(t)$ is convex in $\Del$.
\end{corollary}

\subsection{Representations of $L_w$} \label{sec:proof}

Here and in what follows we use the notations of
Section~\ref{subsec:intro}.

\begin{lemma}
\label{theo:operator} Let $(\Dd{},f,\Del )$ be a regular
lemniscate region and $w=w_{[0]}$ be an analytic function in
$\Dd{}$. Let $L_w(t)=\|w\|^2_t$ (see (\ref{def-s})); then
$L_w(t)\in C^\infty(\Del)$ and for any  $\nu\geq 0$
\begin{equation}
\begin{split}
\label{equ:H-prime}
L_w^{(2\nu+1)}(t)&=\int\limits_{E_{\Dd{},f}(t)}\re
\overline{w}_{[\nu]}w_{[\nu+1]}\,|dz|,\\
L_w^{(2\nu+2)}(t)&=\int\limits_{E_{\Dd{},f}(t)}|w_{[\nu+1]}|^2\;|dz|,
\end{split}
\end{equation}
where $w_{[k]}(z)$ are defined by (\ref{equ:w-def}).
\end{lemma}

\begin{proof}
To apply Lemma~\ref{theo:exhaust} (see also Remark~\ref{rem:nota})
we note that in our case $M=\Com{}$. Let us identify a complex
number $z=x+iy$ with the point $(x,y)\in \R{2}$ so that the
gradient of a real-valued function $h(x,y)$ takes the form
$h'_x+ih'_y$. Then by Cauchy-Riemann theorem we have for any
analytic function $F(z)$
\begin{equation}\label{Cauchy}
\nabla \re F(z)\equiv (\re F(z))'_x+i(\re F(z))'_y=\conj{F'(z)},
\end{equation}
whence
\begin{equation}\label{equ:log-u}
\conj{\vphantom{S^2}\nabla u(z)}=\conj{\vphantom{S^2}\nabla \re
\ln f(z)}=\frac{f'(z)}{f(z)}:=\frac{1}{g(z)}.
\end{equation}
Here and in the sequel we use the logarithms only for brevity of
the computations of gradients and, of course, one can deduce these
formulae directly. Henceforth, we chose the main branch of
logarithm, $\ln 1=0$.

Moreover, we point out that $g(z)$, being defined by
(\ref{equ:log-u}), is an analytic function in $\Dd{}$ because of
the regularity condition: $f'(z)\ne 0$. As a consequence, the same
property is obviously true for the iterates $w_{[k]}$, $k\geq 0$.


Now, let $t\in\Del$. Applying (\ref{equ:log-u}) we get for
$w=w_{[0]}$
$$
L_w(t) =\int\limits_{E_{\Dd{},f}(t)}\frac{|w(z)|^2}{|\nabla
u(z)|}|\nabla
u(z)|\,|dz|=\int\limits_{E_{\Dd{},f}(t)}|w^2(z)g(z)|\,|\nabla
u(z)|\,|dz|.
$$
Hence, substituting $h(z)=|w^2(z)g(z)|$ in (\ref{equ:ex2}) yields
\begin{equation}\label{equ:x2}
L_w'(t)=\int\limits_{E_{\Dd{},f}(t)}\frac{\scal{\nabla
h(z)}{\nabla u(z)}}{|\nabla u(z)|}\, |dz|.
\end{equation}

Further, we note that $ \ln h\equiv \re \ln w^2g, $ and it follows
from (\ref{Cauchy}) that
\begin{equation}\label{equ:h-nabla}
\conj{\nabla h}(z) = h(z)\frac{d{\ }}{dz}(\ln w^2(z)g(z))=
\frac{|w^2g|}{wg}w_{[1]}=\conj{w}w_{[1]}\frac{|g|}{g}.
\end{equation}

On the other hand, we obtain from (\ref{equ:log-u}) and
(\ref{equ:h-nabla})
\begin{equation}\label{equ:scal}
\scal{\nabla u}{\nabla h} = \re \left(\conj{\vphantom{P'}\nabla
u}\nabla h \right) = \frac{1}{|g|}\re \conj{w}w_{[1]}.
\end{equation}
Substituting (\ref{equ:scal}) and  (\ref{equ:log-u}) in
(\ref{equ:x2}) yields
\begin{equation}\label{equ:deriv1}
L_w'(t)=\int\limits_{E_{\Dd{},f}(t)} \re \conj{w}(z)w_{[1]}(z)\,
|dz|.
\end{equation}

To find the second derivative $L_w''(t)$ we notice that $\ln h(z)$
is a harmonic function. It follows then
$$
0=\Delta \ln h(z)=\frac{\Delta h(z)}{h(z)}-\frac{|\nabla
h(z)|^2}{h^2(z)},
$$
and by (\ref{equ:h-nabla}) we arrive at
$$
\Delta h(z)=\frac{1}{h(z)}|\nabla
h(z)|^2=\frac{|ww_{[1]}|^2}{|w^2g|\phantom{{}^2}}=\frac{|w_{[1]}|^2}{|g|}.
$$

After substituting the above expressions in (\ref{equ:ex3}) and
taking into account that $\nabla u=1/|g|$ we obtain from
(\ref{equ:log-u})
\begin{equation}\label{equ:deriv2}
L_w''(t)=\int\limits_{E_{t,\Dd{}}(f)} |w_{[1]}|^2\; |dz|.
\end{equation}

Now we observe that, by virtue of the regularity condition, $w_1$
and the consequent iterates $w_{\nu}$, $\nu\geq2$, are analytic
functions in $\Dd{}$. Since the integral in $L_w''(t)$ takes the
form (\ref{def-s}), it is clear now that formulae
(\ref{equ:H-prime}) can be obtained by induction from
(\ref{equ:deriv1}) and (\ref{equ:deriv2}), and the lemma follows.
\qed
\end{proof}

The further properties of $\|w\|^2_t$ can be deduced by using the
following formalism. Let $(\Dd{},f,\Del)$ be a regular lemniscate
region. We endow the space of all analytic in $\Dd{}$ functions by
a family of the following scalar products
\begin{equation*}\label{equ:scalar}
\scal{u}{v}_{t}:=\int\limits_{E_{t,\Dd{}}(f)}\re \conj{u(z)}v(z)\,
|dz|
\end{equation*}
and $\|w\|^2_t:=\scal{w}{w}_{t}$. Thus, (\ref{equ:H-prime}) can be
rewritten as
\begin{eqnarray}\label{equ:hilbert1}
&D^{2\nu}\|w\|_{t}^2=\|w_{[\nu]}\|_{t}^2, \\
\label{equ:hilbert2} &D^{2\nu+1}\|w\|_{t}^2=
\scal{w_{[\nu]}}{w_{[\nu+1]}}_{t}
\end{eqnarray}
where $D=\frac{d}{dt}$.

We can polarize the preceding identities in the standard way by
using linearity of $D$. Let $\nu=0$ in (\ref{equ:hilbert2}); then
substituting the sum $v_{[0]}+u_{[0]}$ for $w_{[0]}$ yields
\begin{equation*}
\begin{split}
D\|u_{[0]}+v_{[0]}\|^2_{{t}}&=\scal{u_{[0]}+v_{0]}}{u_{[1]}+v_{[1]}}_{t}=\\
&=\scal{u_{[0]}}{u_{[1]}}_{t}+\scal{v_{[0]}}{v_{[1]}}_{t}+
\scal{u_{[0]}}{v_{[1]}}_{t}+ \scal{v_{[0]}}{u_{[1]}}_{t}=\\
&\qquad \qquad \text{(by (\ref{equ:hilbert1}))}
\\
&= D\|u_{[0]}\|^2_{{t}}+ D\|v_{[0]}\|^2_{{t}}+
\scal{u_{[0]}}{v_{[1]}}_{t}+\scal{u_{[1]}}{v_{[0]}}_{t},
\end{split}
\end{equation*}
and it follows that
\begin{equation}\label{equ:D-agree}
2D\scal{u_{[0]}}{v_{[0]}}_{t}=\scal{u_{[0]}}{v_{[1]}}_{t}+\scal{u_{[1]}}{v_{[0]}}_{t}
\end{equation}
holds for two any  analytic functions in $\Dd{}$.

\begin{corollary}\label{theo:M-sys-posit}
Let $(\Dd{},f,\Del )$ be a regular lemniscate region and
$w=w_{[0]}$ be an analytic function in $\Dd{}$. Then for all
$n\geq0$
\begin{equation}\label{equ:bundle}
D^{n}\|w_{[0]}\|^2_t=\scal{w_{[j]}}{w_{[n-j]}}_t,  \quad \forall
j:0\leq j\leq n,\qquad t\in \Del,
\end{equation}
where $w_{[k]}=G_f^k(w)$.
\end{corollary}

\begin{proof}
First we note that by (\ref{equ:hilbert1})-(\ref{equ:hilbert2}) it
is enough to prove that the scalar product in (\ref{equ:bundle})
is independent of $j$. First, applying $k=0$ to
(\ref{equ:hilbert2}), we obtain
\begin{equation}\label{equ:first-M}
D\scalm{w_{[0]}}{w_{[0]}}=\scalm{w_{[0]}}{w_{[1]}}
\end{equation}
which coincides with (\ref{equ:bundle}) for $n=1$, $j=0,1$.

Next, applying $2D$ to  both sides of (\ref{equ:first-M}) we find
by (\ref{equ:D-agree}) that
\begin{equation}\label{equ:with}
2D^2\scal{w_{[0]}}{w_{[0]}}=2D\scalm{w_{[0]}}{w_{[1]}}=
\scalm{w_{[1]}}{w_{[1]}}+\scalm{w_{[0]}}{w_{[2]}},
\end{equation}
hence taking into account (\ref{equ:hilbert1}), $k=1$, we obtain
$$
\scalm{w_{[1]}}{w_{[1]}}=\scalm{w_{[0]}}{w_{[2]}}.
$$
The last identity shows that (\ref{equ:bundle}) is fulfilled for
$n=2$ and $0\leq j\leq 2$.

At the remaining part of  the proof we apply induction over index
$n$. Namely, we suppose that the statement of our assertion holds
for some $n=m\geq 2$. Then applying (\ref{equ:bundle}) for $n=m-1$
we obtain for $u_{[0]}:=w_{[1]}$
\begin{equation*}
\begin{split}
D^{m-1}\;
D^2\|w_{[0]}\|^2_t&=D^{m-1}\|w_{[1]}\|^2_t=D^{m-1}\|u_{[0]}\|^2_t=
\\
&=\scalm{u_{[0]}}{u_{[m-1]}}=\scalm{u_{[i]}}{u_{[m-i-1]}},
\end{split}
\end{equation*}
for all $0\leq i \leq m-1$. Thus, returning to $w_{[0]}$ we get
$$
D^{m+1}\|w_{[0]}\|^2_t=\scalm{w_{[1]}}{w_{[m]}}=\scalm{w_{[2]}}{w_{[m-1]}}=\scalm{w_{[i+1]}}{w_{[m-i]}}
$$
for the same range of $i$. This proves (\ref{equ:bundle}) for
$n=m+1$ and $1\leq j\leq m$.

On the other hand,  applying again the induction assumption $n=m$
 to (\ref{equ:bundle}) we can write
$$
\scalm{w_{[0]}}{w_{[m]}}=\scalm{w_{[1]}}{w_{[m-1]}}
$$
hence applying $2D$ to both sides of the last relation we arrive
at
$$
\scalm{w_{[1]}}{w_{[m]}}+\scalm{w_{[0]}}{w_{[m+1]}}=\scalm{w_{[2]}}{w_{[m-1]}}+\scalm{w_{[1]}}{w_{[m]}},
$$
or $ \scalm{w_{[0]}}{w_{[m+1]}}=\scalm{w_{[2]}}{w_{[m-1]}}$, which
yields (\ref{equ:bundle}) for the remaining cases $j=0$ and
$j=m+1$. The corollary is proved completely. \qed
\end{proof}

We recall the following well-known definition (see, e.g.,
\cite[Ch.~2]{Akhiezer}).

\begin{definition}\label{def:pos}
A sequence $(s_k)$, $k=0,1,\ldots$ is said to be \textit{positive}
if for any polynomial $P(z)=a_0+a_1x+\ldots+a_mx^m$ with real
coefficients which is non-negative  on $\R{}$ it holds
$$
a_0s_0+a_1s_1+\ldots +a_ms_m\geq 0.
$$
\end{definition}

An equivalent definition of  positivity of $(s_k)$
 is that the quadratic forms
$$
\sum_{i,j=0}^m s_{i+j}\xi_i\xi_j
$$
should be positive (semidefinite) for all $m\geq0$
\cite[p.~133]{Widder}. If all the last forms are strictly positive
the sequence $(s_k)$ is called a \textit{strictly} positive
sequence.

We  recall also the following well-known result of H.~Hamburger
\cite{Hamb} (see also \cite[p.~129]{Widder}).

\begin{HT}\label{theo:Hamburger}
A necessary and sufficient condition that there exist at least one
non-decreasing function $\sigma(x)$ such that
    \begin{equation*}\label{equ:Hamburger}
    \int\limits_{-\infty}^{+\infty}x^k\;d\sigma(x)=s_k, \qquad
    (k=0,1,\ldots)
\end{equation*}
with all the integrals converging, is that the sequence $(s_k)$
should be positive.
\end{HT}

Now Theorem ~\ref{theo:H-main} follows from the property
formulated below.

\begin{corollary}\label{corol:theo}
Let $(\Dd{},f,\Del )$ be a regular lemniscate region and $w$ be an
analytic function in $\Dd{}$. Then for all $t\in\Del $ the
sequence
$$
L_{w}^{(k)}(t)=D^k\|w\|^2_t, \qquad k=0,1,\ldots
$$
forms a positive sequence. Moreover, $(L_{w}^{(k)}(t))_{k\geq 0}$
forms a strictly positive sequence (for all $t\in\Del$) if and
only if the system of iterates $\{w_{[k]}\}_{k\geq 0}$ is linearly
independent.
\end{corollary}

\begin{proof}
It follows from (\ref{equ:bundle}) that the corresponding Hankel
matrices (\ref{HAN}) have the form of Gram matrices. Hence,
$(L_{w}^{(k)}(t))_{k\geq 0}$ is a positive sequence by
characteristic property of Gram matrices. The latter assertion of
the corollary now can be proved as follows. Let the sequence
$L_w^{(k)}(t_0)$ fails strict positivity at some $t_0\in\Del$; it
means that there exists an index $N\geq 0$ and a vector
$\xi\in\R{N+1}$, $\xi\ne0$, such that
$$
0=\sum_{i,j=0}^N L_w^{(i+j)}(t_0)\xi_i\xi_j= \sum_{i,j=0}^N
\scal{w_{[i]}}{w_{[j]}}_{t_0}\xi_i\xi_j= \left\|\sum_{i=0}^N
w_{[i]}\xi_i\right\|^2_{t_0}.
$$
But the uniqueness theorem for analytic functions now yields that
$$
\sum_{i=0}^N w_{[i]}\xi_i\equiv 0
$$
on $\Dd{}$, that is the system $\{w_{[i]}\}_{i=0}^{N}$ is linearly
dependent. In particular, it follows that $L_w^{(k)}(t)$ fails
strict positivity for all $t\in\Del$.

The converse property is verified in the same manner and the
assertion is proved. \qed\end{proof}

\subsection{Polynomial lemniscates}\label{subsec:proof}

\begin{proof}[of Theorem~\ref{th-er1}]
To prove the first statement we fix $\Del_j$ to be a principal
regular interval of $P$ and let $w\equiv 1$. In our previous
notations
$$
L_w(t)=|E_P(t)|=\|w\|_t^2
$$
Since $ \Phi_P(t)=\ln |E_{P}(t)|-\frac{t}{n}$, we can examine only
the latter logarithm. Then  convexity of $\ln |E_{P}(t)|$
immediately follows from (\ref{equ:H-prime}) and Cauchy's
inequality:
$$
L_w''(t)L_w(t)-L_w'^2(t)=\|w_{[1]}\|_t^2\|w\|_t^2
-\scal{w}{w_{[1]}}_t^2\geq 0, \qquad t\in\Del_j.
$$
Let now suppose that for some $\tau\in\Del_j$ the equality
$L_w''(\tau)L_w(\tau)-L_w'^2(\tau)=0$ holds. Then $w_{[1]}(z)=
cw(z)$, $z\in E_P(\tau)$, for a constant $c\in\R{}$. Then the
uniqueness theorem for analytic functions yields that the last
identity holds everywhere in $\Dd{}(\Del_j)$ and it follows from
(\ref{equ:w-def}) that
$$
2gw'+g'w=cw
$$
and applying $w\equiv 1$, we get
$$
c=g'(z)\equiv \left(\frac{P(z)}{P'(z)}\right)'.
$$
But this yields $P(z)=c(z-a)P'(z)$, and, consequently (since $P$
is a monic polynomial), $P(z)=(z-a)^n$  with $c=1/n$. Thus, $P$
must be a trivial polynomial and the first statement of the
theorem is proved.

Continuity of $|E_t(P)|$ can be established as follows. We observe
that by virtue of  (\ref{equ:homot}) the following relation holds
\begin{equation}\label{equ:contin}
 |E_{P}(t)|= e^{\frac{t}{n}}|E_{P_t}(0)|
\end{equation}
for all $t\in (-\infty;+\infty]$. On the other hand, we can apply
a result of Eremenko-Hayman \cite[Lemma ~4]{ErHay} which states
that \textit{the lemniscate length $|E_P(0)|$ is a continuous
function of the coefficients of $P$}. Since the coefficients of
$P_t$ (see the explicit expression in (\ref{equ:transf})) are
continuous functions of $t\in (-\infty;+\infty]$, the required
property follows now from (\ref{equ:contin}).

It remains to prove (\ref{limits}). Here we have $t\in
(T_{\nu-1};+\infty)$ where $T_{\nu-1}$ is the largest finite
critical value of $\ln |P|$. Again, applying the Eremenko-Hayman
lemma and (\ref{equ:homot}) we notice that $\lim_{t\to+\infty}
P_t(z)=z^n$ whence
$$
\lim_{t\to+\infty}|E_t(P)|e^{-t/n}=|E_0(z^n)|=2\pi.
$$
and the theorem is proved. \qed\end{proof}

\section{Applications}
\label{sec:appl}

\subsection{$\Dp$-functions}\label{subsec:examples}
Some problems being initially posed by Piranian in \cite{Piran}
and dealing with  monotonicity and convexity of the length
function for $Q_n(z)=z^n-1$ were studied in papers \cite{Butler},
\cite{Elia}. Below we  obtain explicit formulae for the
length-functions $ |E_f(t)| $ for a special class analytic
functions $f$ which include $Q_n$ as a partial case. Our method
involves representation (\ref{equ:s-k-derivative}) to reduce the
problem to a certain hypergeometric differential equation. We
demonstrate it by the following example.

Let $w$ be an analytic function which satisfies the relation
\begin{equation}\label{equ:appl-def}
(\alpha w+\beta w_{[1]})f^\nu=\gamma w+\delta w_{[1]},
\end{equation}
Here we write as above
$
w_{[1]}=G_f( w)=2gw'+g'w.
$
We exclude the trivial case by assuming that
\begin{equation}\label{equ:determinant}
\det \left(%
\begin{array}{cc}
  \alpha & \beta \\
  \gamma & \delta \\
\end{array}%
\right)\ne 0.
\end{equation}

Let $(\Dd{},f,\Del)$ be a regular lemniscate region. To ensure
existence we can suppose that $ \Zer{ f}\ne\emptyset$. Then for
all $t\in\Del$ we have from (\ref{equ:appl-def})
\begin{equation*}\label{equ:appl-diff}
e^{2\nu t}|\alpha w(z)+\beta w_{[1]}(z)|^2= |\gamma w(z)+\delta
w_{[1]}(z)|^2
\end{equation*}
that after integration over $E_{\Dd{},f}(t)$ and using
(\ref{equ:H-prime}) yields
\begin{equation}\label{equ:appl-diff-eq}
(\beta^2 e^{2\nu t}-\delta^2)L_w''(t)+2(\alpha\beta e^{2\nu
t}-\gamma\delta)L_w'(t)+(\alpha^2 e^{2\nu t}-\gamma^2)L_w(t)=0
\end{equation}
with $L_w(t):=\|w\|^2_t$.

Let now choose $w\equiv 1$, such that $L_w(t)=|E_{\Dd{},f}(t)|$
becomes the length function. Then (\ref{equ:appl-def}) can be
rewritten as
$$
f^\nu=\frac{\gamma +\delta g'(z)}{\alpha +\beta g'(z)}.
$$
 Let $z_1\in\Zer{ f}$. We can assume also that $z_1$ is a simple zero
 of $f(z)$, since the arguments similar to that given below show
 that the general case also leads us to the same form of the main equation
 (\ref{equ:appl-g}). Then it follows from the definition of $g$ that
 $g'(z_1)=1$ and consequently we have from (\ref{equ:appl-def}):
$ \gamma+\delta=0, $ whence $\gamma =-\delta$. Moreover, it
follows from (\ref{equ:determinant}) that $\delta\ne0$ and
changing the notations $a=\alpha/\delta$ and $b=\beta/\delta$ we
arrive at
$$
g'(z)=\frac{af^\nu+1}{1-b f^\nu}.
$$
Taking into account that
$$
d(\ln(f/f'))=\frac{dg}{g}=\frac{af^\nu+1}{1-b
f^\nu}\cdot\frac{df}{f}=\frac{af^\nu+1}{1-b f^\nu}\cdot
\frac{df^\nu}{\nu f^\nu}
$$
we obtain after integration and a suitable changing the variables
$\varphi (z)=c_1f(c_2z)$ that
\begin{equation}\label{equ:appl-g}
\varphi '=C(1-\varphi ^\nu)^{\frac{k+1}{\nu}},
\end{equation}
with $k=a/b$.

To analyze the last equation we suppose that $\nu$ is a positive
integer and $k+1\leq \nu$. Then the set of solutions to
(\ref{equ:appl-g}) is still large and such elementary functions as
$\tanh z$, $\sin z$, $z^n-1$, $e^{-z}+1$ satisfy these conditions
(see Table~\ref{tab:1}).

\def\verr{\vphantom{\int\limits_{0}^{0}}}
\def\verrq{\vphantom{\int\limits_{0}^{0}}}

\begin{table}[h]
\begin{center}\begin{tabular}{|r||p{1cm}|p{1cm}|p{2cm}|p{1cm}|}\hline
 $\verr\varphi(z)$  & $k$ &  $\nu$ & $p=\frac{k+1}{2\nu}$ & $C$
   \\\hline\hline
 $\verr 1-z^n$ &    $-\frac{1}{n}$   &   $1$ &
$\frac{n-1}{2n}$ & $-n$
\\\hline
$\verrq\sin z$ &    $0$     &     $2$   & $\frac{1}{4}$ & $1$
\\\hline
$\verrq\tanh z$   & $1$   & $2$ & $\frac{1}{2}$ & $1$ \\\hline
$\verrq e^{-z}+1$   & $1$   & $1$ & $\frac{1}{2}$ & $1$ \\\hline
\end{tabular}
\medskip
\caption{Elementary $\Dp$-functions}\label{tab:1}
\end{center}\end{table}

First, we notice that under our assumptions the integral
$$
F(z)=\int_{0}^{z}\frac{d\z}{(1-\z ^\nu)^{\frac{k+1}{\nu}}},
$$
where the principal branch of the root is chosen, defines a
univalent function $F(z)$ in the unit disk $\Uu=\{z:|z|<1\}$ since
the real part of the derivative
$$
\re F'(z)=\frac{\re (1-z ^\nu)^{\frac{k+1}{\nu}}}{|1-z
^\nu|^{\frac{2(k+1)}{\nu}}}>0, \qquad z\in\Uu
$$
is positive and we can apply the Noshiro-Warschawski theorem
\cite[p.~47]{Duren}. Thus, for any $C\ne 0$ the function
$\varphi(z)$ given by
$$
F(\varphi(z))=Cz
$$
is analytic (and univalent) in $S:=\frac{1}{C}F\left(\Uu\right)$,
$\varphi:S\to \Uu $ and satisfies (\ref{equ:appl-g}) there. We
call such a function $\varphi$ a $\Dp$-\textit{function}.

An important property of $\Dp$-functions is that their critical
values have the same magnitude: $\varphi'(\z)=0\Rightarrow
|\varphi(\z)|=1$. Clearly, $(S^*,\varphi,(-\infty,0))$ is a
regular lemniscate region of the function $\varphi$, where
$S^*=S\setminus \{0\}$.

Another representation of $F$ can be easily found by using the
Gauss hypergeometric function
$$
F(\z)=\zeta \;{}_2F_{1} \left(
\frac{1+k}{\nu},\frac{1}{\nu};\frac{1+\nu}{\nu},{\zeta}^{\nu}
 \right), \qquad \z\in\Uu.
$$

\begin{theorem}\label{theo:butler}
Let $\nu\in\Z{+}$ and $k\leq \nu-1 $. Then in the preceding
notations, the following formula holds
\begin{equation}\label{equ:appl-exact1}
|E_{S^*,\varphi}(t)|= \frac{2\pi e^t}{|C|} \;{}_2F_1(p,p;1;e^{2\nu
t}),
\end{equation}
where $p=(k+1)/2\nu$ and ${}_2F_1(a,b;c;z)$ is the Gauss
hypergeometric function.
\end{theorem}

\begin{proof}
Let $\varphi (z)$  a $\Dp{}$-function with parameters $\nu$, $k$,
and $C$. Then
$$
g'_\varphi (z)\equiv
\frac{\varphi(z)}{\varphi'(z)}=\frac{1+k\varphi ^\nu}{1-\varphi
^\nu},
$$
and after comparing with (\ref{equ:appl-def}) we have $ \alpha=k$,
$\beta=1$, $\gamma=-\delta=-1$. We notice that in  our previous
notations the function $ |E_{S^*,\varphi}(t)|=L_w(t)$ is a
solution of (\ref{equ:appl-diff-eq}) where $w\equiv 1$. Hence,
letting $t=\ln\tau$
 and
 $$
 L(\tau):=|E_{S^*,\varphi}(\ln \tau)|
 $$
 one can readily see from
(\ref{equ:appl-diff-eq}) that
\begin{equation}\label{equ:appl-L}
\tau^2(\tau^{2\nu}-1)L''(\tau)+[1+(2k+1)\tau^{2\nu}]\tau
L'(\tau)+(k^2 \tau^{2\nu}-1)L(\tau)=0.
\end{equation}

Further, the change of variables $x=\tau^{2\nu}$ and $
L(\tau)=\tau y(\tau^{2\nu})$ reduces (\ref{equ:appl-L}) to the
hypergeometric canonic form
\begin{equation}\label{equ:appl-canonic}
x(1-x)y''(x)+\left[1-x(1+\frac{k+1}{\nu})\right]
y'(x)-\left(\frac{k+1}{2\nu}\right)^2y(x)=0.
\end{equation}
The general solution to (\ref{equ:appl-canonic}) in $[0;1)$ can be
written  as follows
$$
y(\tau)= {}_2F_1(p,p;1;\tau)\;\lambda+\
{}_2F_1(p,p;2p;1-\tau)\;\mu{},
$$
where $\lambda,\mu\in \R{}$ (see \cite[\S~2.3.1]{TransFunct}).

Taking into account that $z=0\in S$ is a simple zero of $\varphi$
we conclude that the length function
$|E_{S^*,\varphi}(t)|=O(e^{t})$ as $t\to-\infty$. Thus, $y(\tau)$
is a bounded function near $\tau=+0$. Since
$$
\lim_{\tau\to +0}{}_2F_1(p,p;2p;1-\tau)=\infty
$$
(see \cite[\S~2.1.3]{TransFunct})) we have $\mu=0$ and after
substitution of the old notations we arrive at
$$
|E_{S^*,\varphi}(t)|= \lambda e^t\;{}_2F_1(p,p;1; e^{2\nu t}).
$$
The precise form of $\lambda$ can now be found by the asymptotic
behavior
$$
|E_{S^*,\varphi}(t)|=\frac{2\pi e^t }{|\varphi'(0)|}, \qquad t\to
-\infty,
$$
On the other hand, by the definition of $\Dp$-function
$|\varphi'(0)|=|C|$ and the theorem is proved. \qed\end{proof}

\begin{remark}
In the  case $t>0$, the last theorem is still meaningful provided
that the level sets $\{|\varphi(z)|=e^t\}$ are compact or
\textit{periodic} curves (the last, e.g., corresponds to
$\varphi=\sin z$). These cases are treated in \cite{KT}.
\end{remark}

\subsection{The structure of $\sigma^{P,\Del{}}$} \label{sec:measure}
Here we describe an explicit structure of the measure function
(for monic polynomials) $\sigma^{P,\Del_j}(x)$  when $j=\nu$, i.e.
the interval $(T_{\nu-1},+\infty)$ is free of critical values of
$\ln |P|$. The latter means that the corresponding lemniscates
$E_P(t)$, $t\in\Del_\nu$, are single-component closed curves.


\begin{theorem}\label{theo:shimorin}
Let $P$ be a monic polynomial of degree $n$ and
$$
T(P):=T_\nu=\max_{P'(\z_k)=0} \ln|P(\z_k)| $$ is the largest
singular value. Then for all $t\geq T(P)$ the following
representation holds
\begin{equation}\label{equ:shim-rep}
|E_P(t)|=2\pi e^{t/n}\left(1+
\sum_{k=2}^{+\infty}|c_{-k}|^2e^{-2kt/n}\right),
\end{equation}
 and
\begin{equation}\label{qq}
\sum_{k\geq 2}|c_{-k}|^2e^{-2kT(P)/n}<+\infty.
\end{equation}
Here $c_k$ is the $k$th Laurent coefficient of
$\sqrt{\varphi'(\z)}$ near the infinity, where
$P(\varphi(\z))=\z^n$.
\end{theorem}

\begin{proof}
Let $t>T(P)$ be chosen arbitrary. Then $E_{P}(t)$ is a simple
Jordan curve which is the boundary of a simply-connected domain
$$
D_t=\{z\in \Com{}: |P(z)|<e^t\}.
$$ Let
$D^*_t=\overline{\Com{}}\setminus \overline{D_t}$. Then $P(z)$
maps $D^*(t)$ onto $\Uv _t:=\{\z\in\overline{\Com{}}: |\z|>e^t\}$.
Since $P'(z)\ne 0$ in $\Uv _t$ the analytic function
$F(z)=P^{1/n}(z)$, $F(z)\sim z$ as $z\to\infty$, is well defined
in $D^*(t)$. Moreover, $F(z)$ is univalent in $D^*(t)$ for all
$t\geq T(P)$. We denote by $\varphi(\z)$ the inverse function
which is also a univalent function and observe that
\begin{equation}\label{equ:shim2}
\varphi:\Uv _{T(P)/n}\to D^*(T(P)); \qquad \varphi'(\z)\ne 0 ,
\quad \z\in \Uv _{T(P)/n}.
\end{equation}
We have $\varphi(\z)\sim \z$ as $\z\to\infty$ and letting
$P(z)=z^n+a_{1}z^{n-1}\ldots+a_n$ we obtain
\begin{equation}\label{equ:Sigma}
\varphi'(\z)=1-\frac{(n-1)a_1^2-2na_2}{2n^2}\frac{1}{\z^2}+\ldots,
\qquad \z\to \infty.
\end{equation}
It follows that $\sqrt{\varphi'}$ is a well-defined analytic
function in $\overline{\Uv }_{T(P)/n}$ (for sake of completeness
let $\sqrt{1}=1$). Thus, it can be expanded in the Laurent series
\begin{equation}\label{equ:Laurent}
\sqrt{\varphi'(\z)}=1+\sum_{k=2}^{+\infty}\frac{c_{-k}}{\z^k},
\qquad |\z|>e^{T(P)/n}.
\end{equation}

Next, we observe that for $t\geq T(P)$, the curve $E_{P}(t)$ is
homeomorphic to a circle and can be naturally parameterized by
$$
E_{P}(t)=\{\varphi(\theta): \quad \theta \in e^{t/n}\Tcirc\},
$$
where $\Tcirc=\partial \Uu$ is the unit circle,  and we have for
the length function
\begin{equation}\label{equ:len1}
\begin{split}
&|E_{P}(t)|=\int\limits_{e^{t/n}\Tcirc}|\varphi'(\z)|\,|d\z|
=\int\limits_{e^{t/n}\Tcirc}|\sqrt{\varphi'(\z)}|^2\,|d\z|=\\
&=\int\limits_{e^{t/n}\Tcirc}\sum_{k=-\infty}^{+\infty}|c_k|^2|\z|^{2k}\,|d\z|
= 2\pi e^{t/n}\left(1+
\sum_{k=2}^{+\infty}|c_{-k}|^2e^{-2kt/n}\right)
\end{split}
\end{equation}
where convergence of the corresponding series for $t>T(P)$ follows
from that one in (\ref{equ:Laurent}).

To establish (\ref{qq}) we note that it follows from
(\ref{equ:len1}) that $|E_{P}(t)|$ is a decreasing function in
$(T(P);+\infty)$ and continuous in $[T(P),+\infty)$ (see
Theorem~\ref{th-er1}). Thus, the right-hand side of
(\ref{equ:len1}) is at most $|E_{P}(T(P))|$ for all $t>T(P)$. Then
it follows from positivity of the terms of the corresponding
series and the mentioned continuity that (\ref{equ:len1}) is still
true for $t=T(P)$, which in turn yields (\ref{qq}) and completes
the proof. \qed\end{proof}

\begin{corollary}\label{cor:low-est}
In the notations of Theorem~\ref{theo:shimorin} we have the
following lower estimate
\begin{equation}\label{equ:low-est}
e^{4t/n}\left(\frac{e ^{-t/n}}{2\pi} |E_{P}(t)|-1\right)\geq
\left|\frac{(n-1)a_1^2-2na_2}{4n^2}\right|^2 , \qquad t\geq T(P)
\end{equation}
where $P(z)=z^n+a_{1}z^{n-1}\ldots+a_n$ and the estimate is sharp.
\end{corollary}

\begin{proof}
The estimate easily follows from (\ref{equ:Sigma}). On the other
hand, the same formula shows that the left-hand side of
(\ref{equ:low-est}), denote it by $h(t)$, is a decreasing function
as $t\geq T(P)$ and
$$
\lim_{t\to+\infty}
h(t)=|c_{-2}|^2=\left|\frac{(n-1)a_1^2-2na_2}{4n^2}\right|^2
$$
which proves a sharp character of (\ref{equ:low-est}).
\qed\end{proof}

It is helpful also to notice that the expression in the right-hand
side of (\ref{equ:low-est}) has the form
$$
\frac{(n-1)a_1^2-2na_2}{4n^2}=\frac{1}{4n^2}\sum_{1\leq i<j\leq
n}(z_i-z_j)^2
$$
where $P(z)=\prod_{j=1}^n(z-z_j)$.

 We recall  (see \cite[p.~145]{Widder}) that a function $f(t)$ is
\textit{completely monotonic} in $(a,b)$ if it has non-negative
derivatives of all orders there:
\begin{equation}\label{equ:cm}
(-1)^kf^{(k)}(t)\geq 0.
\end{equation}
A function $f(t)$ is said to be completely monotonic in $[a,b)$ if
it is continuous there
 and satisfies (\ref{equ:cm}) in $(a,b)$.

\begin{corollary}\label{cool:cm}
$|E_{P}(t)|e^{-t/n}$ is completely monotonic in $[T(P),+\infty)$.
\end{corollary}

Finally, we briefly discuss the mentioned in
Corollary~\ref{corol:cont} analytic continuation property. In the
post critical case $t>T(P)$ this fact can be established directly.
Indeed, let
\begin{equation}\label{equ:analytic}
F(z)=2\pi e^{z}\left(1+
\sum_{k=2}^{+\infty}|c_{-k}|^2e^{-2kz}\right), \qquad \re z\geq
\frac{T(P)}{n},
\end{equation}
then it follows from (\ref{qq}) that $F(\z)$ is a single-valued
analytic function and
\begin{equation}\label{complexific}
|E_{P}(t)|=F\left(\frac{t}{n}\right), \qquad t\in \R{}
\end{equation}
is a desirable continuation.

Moreover, we note that $F(z)$ is a $2\pi i$-periodic function. Let
$$
\lambda(\z)=F(\ln \z)=2\pi \z\left(1+
\sum_{k=2}^{+\infty}\frac{|c_{-k}|^2}{\z^{2k}}\right), \qquad
|z|>e^{\frac{T(P)}{n}}
$$
Then the last formula shows that $\lambda(\z)$ is an \textit{odd}
analytic function and
\begin{equation}\label{equ:lam-e}
|E_{P}(t)|=\lambda(e^{t/n}), \quad \forall t\geq T(P).
\end{equation}

\begin{acknowledgement}  The authors are grateful to Alexander
Eremenko, Bj\"orn Gus\-taf\-sson, Henrik Shahgholian, Harold
Shapiro and Serguei Shimorin for their constructive comments. We
are also indebted to Klaus Steffen and the referee for careful
reading and suggestions which led to a considerable improvement of
the paper.
\end{acknowledgement}


\begin{thebibliography}{AAAA}
\bibitem{Akhiezer}
Akhiezer, N.I.: \ The Classical Moment Problem and Some Related
Questions in Analysis, English translation. Oliver and Boyd,
Edingburgh, 1965

\bibitem{Alzer}
Alzer, H., Berg, C.: Some classes of completely monotonic
functions. \textit{Ann. Acad. Sci. Fen. Math.},
\textbf{27},445--460 (2002)

\bibitem{semi}
Berg, C., Christensen, J.P.R., Ressel, P.: \ Harmonic Analysis on
Semigroups. New-York-Berlin, Springer-Verlag. 1984


\bibitem{BergDuran}
Berg, C., Duran, A.J.: A transformation from Hausdorff to
Stieltjes moment sequences (\textrm{to appear})


\bibitem{Bern-AMF}
Bernstein, S.N.: Sur les fonctions absolument monotones.
\textit{Acta math}. \textbf{52}, 1--66 (1928)

\bibitem{Bor}
Borwein, P.: The arc length of the lemniscate $\{|P(z)|=1\}$.
\textit{Proc. Amer. Math.Soc.}, \textbf{123}, 797--799 (1995)


\bibitem{Butler}
Butler, J.P.: The perimeter of a rose. \textit{Amer. Math.
Monthly.} \textbf{98}, no. 2, 139--143 (1991)

\bibitem{Duren}
Duren, P.L.: Univalent functions, Grundlehren der Mathematischen
Wissenschaften, vol. 259, Springer-Verlag, New York (1983).

\bibitem{Elia}
Elia, M., Galizia Angeli, M.T.: The length of a lemniscate.
\textit{Publ. Inst. Math. Beograd (N.S.)} \textbf{36}, 51--55
(1984)
\bibitem{TransFunct}
Erdelyi, A. (ed.): Higher transcendental functions. Vol. I.
Bateman Manuscript Project, California Institute of Technology.
Malabar, Florida: Robert E. Krieger Publishing Company. XXVI,
1981.

\bibitem{Erdelyi}
Erd\'elyi, T.: Paul Erd\"os and polynomials. \textit{Jour. of
Appr. Theory}, \textbf{94}, 2--14 (1998)

\bibitem{Erd}
Erd\"os, P.: Some old and new problems in approximation theory:
research problem 95-1. \textit{Constr. Approx.} \textbf{11},
419-421 (1995).

\bibitem{EHP}
Erd\"os, P., Herzog, F., Piranian, G.: Metric properties of
polynomials. \textit{J. D'Analyse Math.}, \textbf{6}, 125--148
(1958)



\bibitem{ErHay}
Eremenko, A., Hayman, W.: On the length of lemniscates.
\textit{Mich. Math.~ J.}, \textbf{46}, 409--415 (1999)

\bibitem{Federer}
Federer, H.: Geometric Measure Theory. Classics in Mathematics.
Berlin: Springer-Verlag. xvi, 1996

\bibitem{GPP}
Giordano, C., Palumbo, B., Pe\v{c}ari\`{c}, J.: Remarks on the
Hankel determinants inequalities. \textit{Rend. Circ. Mat. Di
Palermo}, Serie II, \textbf{XLVI}, 279-286 (1997).

\bibitem{Hamb}
Hamburger, H.: \"Uber eine Erweiterung des Stieltjesschen
Momentenproblems. \textit{Math. Ann.}, \textbf{81}, (1920).



\bibitem{Hille}
Hille, E.: Analytic function theory. Vol.~II. Ginn \& Co. New
York, 1962.

\bibitem{Kimb}
Kimberling, C.H.: A probabilistic interpretation of complete
monotonicity. \textit{Aequations Math.}, \textbf{10}, 152-164
(1974).

\bibitem{Klya}
Klyachin, V. A.: New examples of tubular minimal surfaces of
arbitrary codimension. \textit{Math. Notes} \textbf{62}, no. 1-2,
129--131 (1998)


\bibitem{KT}
Kuznetsova, O.S., Tkachev, V.G.: Analysis on lemniscates and
Hamburger's moments. Preprint. TRITA-MAT-2003-04. Division of
Math., Royal Inst. of Techn., Stockholm, (2003)


\bibitem{Marden}
Marden, M.: The Geometry of Zeros of a Polynomial in a Complex
Variable. Vol. 3 of Mathematics Surveys. Amer. Math. Soc.,
Providence, 1949



\bibitem{Mik79}
Mikljukov, V.M.: Some properties of tubular minimal surfaces in $
R\sp{n}$. \textit{Dokl. Akad. Nauk SSSR} \textbf{247}, no. 3,
549--552 (1979).

\bibitem{MT89}
Miklyukov, V.M., Tkachev, V.G.: Some properties of tubular minimal
surfaces of arbitrary codimension. \textit{Math. USSR, Sb.}
\textbf{68}, no.~1, 133--150 (1990)


\bibitem{Samko}
Miller, K.S, Samko, S.G.: Completely monotonic functions.
\textit{Integr. transf. and special funct.,} \textbf{12}, no~4,
389-402 (2001)


\bibitem{Piran}
Piranian, G.: The length of a lemniscate. \textit{Amer. Math.
Month.}, \textbf{87}, 555--556 (1980)

\bibitem{Pom59-2}
Pommerenke, Ch.: On some problems of Erd\"os, Herzog and Piranian.
\textit{Mich. Math. J.} \textbf{6}, 221--225 (1959)

\bibitem{Pom61}
Pommerenke, Ch.: On some metric properties of polynomials II.
\textit{Mich. Math. J.} \textbf{8}, 49--54 (1961)



\bibitem{Umem}
Umemura, Y.: Measures on infinite dimensional vector spaces.
\textit{Publ RIMS. Kyoto Univ.} \textbf{1}, 1--47 (1965)



\bibitem{Widder}
Widder D.V., The Laplace Transform. Princeton, University Press,
Princeton (1946)


\end{thebibliography}
\end{document}